\documentclass[12pt]{article}
                                                                               
 \usepackage{amsfonts,amssymb,eucal}

\usepackage{amsthm}
 \usepackage{amsmath}
\usepackage{amscd}
 \usepackage{latexsym} 
\numberwithin{equation}{section}

 \newcommand{\mb}[1]{{\mbox{\boldmath{$#1$}}}}% mathematical bold
  \newcommand{\mc}[1]{{\mathcal{#1}}}% mathematical caligraphic
 \newcommand{\got}[1]{{\mathfrak{#1}}}% gothic with mbox for  mathematic
\newcommand{\db}[1]{{\mathbb{#1}}}% double
\newcommand{\gata}{\square} 
\newcommand{\pa}{\partial}
\newcommand{\C}{\ensuremath{\mathbb{C}}}
\newcommand{\N}{\ensuremath{\mathbb{N}}}

 \newcommand{\LL}{\ensuremath{\mathbb{L}}}
\newcommand{\T}{\ensuremath{S(U(1)\times U(1) \times U(1))}}
 \newcommand{\Hi}{\ensuremath{\mathcal{H}}}% Hilbert space
 
 \newcommand{\Hinf}{\ensuremath{\mathcal{H}^{\infty}}}% \infty Hilbert space
  \newcommand{\U}{\ensuremath{\mathcal{U}}}% Universal
  \newcommand{\g}{\ensuremath{\got{g}}}% Lie algebra g
\newcommand{\m}{\ensuremath{\got{m}}}% Lie algebra m
\newcommand{\bb}{\ensuremath{\got{b}}}% Lie algebra b
\newcommand{\gc}{\ensuremath{\got{g}_{\C}}}% Lie algebra g-complexificat
\newcommand{\Ug}{\ensuremath{\mathcal{U}({\g}})}% Universal algebra of
\newcommand{\Ugc}{\ensuremath{\mathcal{U}({\gc}})}% Universal algebra of
\newcommand{\Ub}{\ensuremath{\mathcal{U}({\bb}})}% Universal algebra of
\renewcommand{\P}{\ensuremath{\mathbb{P}}}
 \newcommand{\Ph}{\ensuremath{\P (\Hi )}}% Projective of Hilbert
\newcommand{\Phinf}{\ensuremath{\P (\Hinf )}}% Projective of Hilbert  % infty
\newcommand{\DM}{\ensuremath{{\got{D}}_M }}% D_M- module
\newcommand{\AM}{\ensuremath{{\got{A}}_M }}% A_M- module
 \newcommand{\D}{\ensuremath{{\got{D}}}}% D-module% now called sheaf
\newcommand{\A}{\ensuremath{{\got{A}}}}% A-module
\newcommand{\AAA}{\ensuremath{{\db{A}}_M}}% A-module
\newcommand{\am}{\ensuremath{{\bf{A}}_M}}% 
 
\newcommand{\FSB}{symmetric Fock space }
 
 \newtheorem{Remark}{Remark}

\newcommand{\fl}{\ensuremath{{\mathcal{F}}_{\Hi}}}%Bargmann H space
\newtheorem{Proposition}{Proposition}
 
\theoremstyle{definition}%nou 

\newtheorem{deff}{Definition}%martin

\begin{document}
%{\Huge$${\got{The~ National~ Conference~}}$$
%$${\got{ on}}$$
%$${\got{ Theoretical~ Physics}}$$}\\[2cm]
%\begin{center}
%{\LARGE {\bf September 13-16, 2002\\[3ex] Bucharest}}\\[3cm]
\begin{center}
{\Large {\bf Differential operators on orbits of coherent states}}\\[2ex]
S. Berceanu, A. Gheorghe\\[2ex]
 National 
 Institute for Physics and Nuclear Engineering\\
         Department of Theoretical Physics\\
         PO BOX MG-6, Bucharest-Magurele, Romania\\
         E-mail: Berceanu@theor1.theory.nipne.ro; Cezar@theor1.theory.nipne.ro\\
\end{center}
\begin{abstract}
We emphasize some properties of coherent state groups, i.e. groups
whose quotient with the stationary groups, are manifolds which admit a 
holomorphic embedding in a projective Hilbert space. We determine the 
differential action of the generators of the representation of coherent
state groups on the symmetric Fock space attached to the dual of the
Hilbert space of the  representation. This permits a realization of
coherent state Lie algebras  by first-order differential operators
with holomorphic polynomial coefficients on  K\"ahler coherent state orbits. 
\end{abstract}

\section{Introduction}

 The differential action of the
 generators of the groups  on  coherent state manifolds
  which have the structure of
  hermitian symmetric spaces can be written down as a sum
 of two terms, one a polynomial $P$, and the second one  a sum of partial
 derivatives times some  polynomials $Q$-s,
    the degree of  polynomials  being  less than
 3 \cite{sbcag,sbl}.
It is interesting to investigate the same problem as in
 \cite{sbcag,sbl} on flag manifolds 
  \cite{bfr}.  Some
  results are available \cite{dob}, but they are not easily  handled. 

 Our investigations on the differential action of the
generators of hermitian  groups on hermitian symmetric spaces
have been extended 
to semisimple Lie groups acting   on coherent state manifolds
which admit a K\"ahler structure, and  
 explicit formulas for the polynomials $P$ and $Q$-s  have been given
 \cite{sbcpol}.  Explicit formulas for the
simplest example of a compact nonsymmetric coherent state     
 manifold, $SU(3)/\T$, where the degree of the
polynomial is already 3, have been also obtained  \cite{sbcpol}. Here we
discuss in the context of the so called
coherent state (shortly, CS)-groups \cite{lis1,lis2,lis,neeb}
 the space of functions on which these differential operators
act. 

We emphasize some properties of CS-groups,
i.e. groups whose quotient with the stationary groups are manifolds
which admit a holomorphic embedding in a projective Hilbert
space. This class of groups contains all compact groups, all simple
hermitian groups, certain solvable groups and also mixed groups as the
semidirect product of the Heisenberg group and the symplectic group
\cite{neeb}.   
We determine the differential action of the generators of the
representation of the CS-group  on the symmetric Fock
space attached to the dual of the Hilbert space of
the representation.

  The 
 coherent states are a  useful tool of investigation of quantum and classical
 systems \cite{perG}. It was shown in \cite{sbcag,sbl} that a linear 
Hamiltonian in the
 generators of the groups implies  equivalent quantum and classical
evolution. It  was proved that for Hermitian symmetric spaces the evolution
 equation generated by Hamiltonians which are linear in the generators of the
 group is a matrix Riccati equation.
 It is interesting to see how it looks like
 the corresponding equation of
 motion generated by linear Hamiltonians for CS-manifolds. So, the present
 paper gives the technical tools for this calculation. 
The degrees of the polynomials $P$ and $Q$-s
 for flag manifolds are grater
 than 2. In references \cite{sbctim} we have underlined
the realization of coherent state algebras by differential operators,
giving explicit formulas for the semisimple case.

 Another field of possible applications is the determination of the Berry
   phase \cite{berry} on CS-manifolds.
 In \cite{sbcag,sbl} there were presented explicit expressions for the 
 Berry phase for the complex Grassmann manifold. These results 
were used further 
 \cite{sb7} for explicit calculation of the symplectic area of geodesic
 triangles on the complex Grassmann manifold and its noncompact dual.
Also we have considered explicit boson expansions
 for collective models on K\"ahler CS-orbits in
\cite{sbcbuc}.

The general framework in which this paper must be considered is the deep
 relationship between coherent states and geometry \cite{sbcl}. 
 
 The paper is laid out as follows. \S \ref{CSR} contains
the definition of  CS-orbits, 
 in the context of Lisiecki \cite{lis1,lis2,lis} 
and Neeb \cite{neeb}. The geometry of coherent state manifolds for
compact groups was previously considered in \cite{morse}.
\S \ref{CSvectors} deals with the so called
Perelomov's CS-vectors.  The coherent vectors are  defined
taking into account that CS-representations are realized by highest
weight representations and the manifold of coherent states is a
reductive
 homogeneous  space. In \S \ref{fock} we construct the space of
functions on which the differential operators will act.
In \S \ref{DIFF} we study the representations of Lie algebras of
CS-groups by differential operators. 

We use for the scalar product the convention:
$(\lambda x,y)=\bar{\lambda}(x,y)$, $ x, y\in\Hi ,\lambda\in\C $.
\section{CS-representations}\label{CSR}

Let us consider the triplet $(G, T, \Hi )$, where $T$ is
 a continuous, unitary
representation 
 of the  Lie group $G$
 on the   separable  complex  Hilbert space \Hi .
Let us denote by $\Hinf$ the dense subspace of \Hi~ consisting of
 those vectors
$v$ for which the orbit map $G\rightarrow \Hi , g\mapsto
 T (g).v$ is smooth.  Let us
pick up $e_0\in \Hinf$ and let  the notation:
%\begin{equation}\label{unu1}
$e_{g,0}:=T(g).e_0, g\in G$.
%%%%%\end{equation}
We have an action $G\times \Hinf\rightarrow\Hinf$, $g.e_0 :=
e_{g,0}$. When there is no possibility of confusion, we write just
$e_{g}$ for $e_{g,0}$. 
Let us denote by  $[~]:\Hi^*:=\Hi\setminus\{0\}\rightarrow\Ph=\Hi^*/ \sim
$ the projection with respect to the equivalence relation
 $[\lambda x]\sim [x],~ \lambda\in \C^*,~x\in\Hi^*$. So,
$[.]:\Hi^*\rightarrow \Ph , ~[v]=\C v$. The action 
$ G\times \Hinf\rightarrow\Hinf$ extends to the action
$G\times \Phinf\rightarrow\Phinf , g.[v]:=[g.v]$. 

Let us now denote by $H$  the isotropy group $H:=G_{[e_0]}:=
\{g\in G|g.e_0\in\C e_0\}$.
We shall consider (generalized) coherent 
 states on complex  homogeneous manifolds $M\cong
G/H$, imposing the restriction that $M$ be  a complex submanifold of
\Phinf .

\newpage
\begin{deff}\label{def1}
 a) The orbit $M$ is called a CS-{\em orbit} if 
there exits  a holomorphic embedding
$\iota : M \hookrightarrow \Phinf$.
In such a case $M$ is also called CS-{\it manifold}.

b) $(T,\Hi )$ is called a  CS-{\em representation} if there exists a cyclic
 vector $0\neq
e_0\in\Hinf $ such that $M$ is a CS-orbit.

c)  The groups $G$ which admit
CS-representations  are called CS-{\it groups}, and their Lie algebras
\g~ are called CS-{\it Lie algebras}.
\end{deff}

The  $G$-invariant complex structures on the homogeneous
 space $M=G/H$   
 can be introduced in an algebraic manner. 
 For $X\in\g$, where \g ~is the Lie algebra of the Lie group $G$,  let
us define the unbounded operator $dT(X)$ on \Hi~ by
$dT(X).v :=\left. {d}/{dt}\right|_{t=0} T(\exp tX).v$
whenever the limit on the right hand side exists. We obtain a
representation of the Lie algebra \g~ on \Hinf , {\it the derived
representation}, and we denote
${\mb{X}}.v:=dT(X).v$ for $X\in\g ,v\in \Hinf$. Extending $dT$ by complex
linearity, we get a representation of the complex Lie algebra \gc~ on
the complex vector space \Hinf . 
Lemma XV.2.3 p. 651 in \cite{neeb}
and Prop. 4.1 in \cite{lis} determine when a smooth vector
generates a complex orbit in \Phinf .
Let now  denote by $B:=<\exp_{G_{\C}} \bb >$ the Lie group
 corresponding to the Lie algebra
\bb ,  with
$\got{b}:=\overline{\got{b}(e_0)}$, where $\got{b}
 (v):=\{X\in\gc :X.v\in \C v\}=(\g_{\C})_{[v]}$. The group $B$ will be
supposed  to be closed in the complexification $G_{\C}$ of
$G$, and in fact this assumption is justified for CS-groups $G$
(cf. Lemma XII.1.2. p. 495 in \cite{neeb}). 
 Then
the  complex structure on $M$ is induced by an 
embedding in a complex manifold,
 $i_1:M\cong G/H \hookrightarrow  G_{\C}/B$.
 We  consider such
manifolds  which admit a holomorphic
 embedding $i_2: G_{\C}/B\hookrightarrow
\Phinf$.   Then the embedding $\iota
=i_1\circ i_2$,  $\iota : M  \hookrightarrow \Phinf$
 is a holomorphic embedding,  in the
sense that the complex structure comes as in 
  Theorem XV.1.1 and Proposition XV.1.2 p. 646 in \cite{neeb}.

\section{CS-vectors }\label{CSvectors}

Now we construct what we  call  Perelomov's
(generalized)  coherent state vectors, or simply CS-vectors,  based on the
homogeneous manifold $M\cong G/H$. Usually \cite{perG},
this construction is done for (semi)simple Lie groups $G$ with $H:=K$,
where $K$ is a
maximally compact subgroup of $G$. Here we  do this
construction for the  CS-groups $G$ in the meaning of
Definition \ref{def1}.

  We denote also  by $T$ the holomorphic extension of the
representation $T$ of $G$ to  the complexification $G_{\C}$ of $G$,
 whenever this
holomorphic extension exists. In fact, it can be shown that
 in the situations under interest
in this paper, this holomorphic extension exists 
\cite{neeb95,neeb96}.
 Then there exists the 
homomorphisms $\chi_0 $ ($\chi$), $\chi_0: ~H\rightarrow \db{T}$,
($\chi : B\rightarrow \C^*$), 
such that 
$ H = \{g\in G| e_g=\chi_0(g)e_0\}$ (respectively, 
$B = \{g\in G_{\C}| e_g=\chi(g)e_0\})$, where
$\db{T}$ denotes the torus $\db{T}:=\{ z\in\C| |z|=1\}$. 

For the homogeneous space $M=G/H$ of cosets $\{gH\}$, let $\lambda
:G\rightarrow G/H$ be the natural projection $g\mapsto gH$, and let
$o:=\lambda (\mb{1})$, where $\mb{1}$ is the unit element of $G$. 
 Choosing a section $\sigma :G/H\rightarrow
G$ such that $\sigma ( o )={\bf{1}} $, every element $g\in G$ can
be written down as $g=\tilde{g}(g)h(g)$, where
$\tilde{g}(g)\in G/H$ and $h(g)\in H$. Then we have
$e_g=e^{i\alpha (h(g))}e_{\tilde{g}(g)}$,
where $e^{i\alpha(h(g))}=\chi_{0}(h)$. Now we take into account that
$M$ also admits an embedding in  $G_{\C}/B$. We choose a local system
of coordinates parametrized by $z_g$ (denoted also simply $z$, where
there is no possibility of confusion) on  $G_{\C}/B$. Choosing a
section 
$G_{\C}/B
\rightarrow G_{\C}$ such that any element $g\in G_{\C}$ can be written
as $g=\tilde{g}_bb(g)$, where $\tilde{g}_b\in G_{\C}/B$, and
$ b(g)\in B$, we have
$e_g=\Lambda (g)e_{z_g}$,
where
$\Lambda (g)= \chi
(b(g))=e^{i\alpha(h(g))}(e_{z_g},e_{z_g})^{-\frac{1}{2}}$.

 Let    us denote by
$\got{m}$
 the vector space orthogonal to $\got{h}$ 
 of
the Lie algebra \g , i.e. we have the vector space
decomposition $\g=\got{h}+\got{m}$. Even more, it can be shown that the
 vector space decomposition $\g=\got{h}+\got{m}$ is
Ad $H$-invariant. The homogeneous spaces $M\cong G/H$ with this
decomposition are called 
 {\em reductive spaces} (cf. \cite{nomizu})
and it can be proved that the
CS-manifolds  are  reductive spaces. More exactly,
using
  Lemma
III.2.(iii) in \cite{neeb1} and Lemma
XV.2.5 p. 652 in \cite{neeb},  it
can be proved that:
\begin{Remark}\label{red}  The homogeneous
 coherent state manifold $M\cong G/H$, for
which the isotropy representation has discrete kernel, or for
admissible Lie algebras and faithful CS-representations, is a
reductive space.
\end{Remark}
So, {\it the tangent space to $M$ at $o$
can be identified with $\got{m}$}.
Now  remember (cf. Proposition XV.2.4 p 651, Proposition
XV.2.6 p. 652  in \cite{neeb} and  Theorem XV.2.10 p. 655
in \cite{neeb}, where  the
 algebra \g~ is
supposed to be admissible  (cf. Definition VII.3.2
at p. 252)) that for CS-groups, the
 CS-representations are highest weight representations
(cf. Definition X.2.9 p. 399 in \cite{neeb}),
  and the vector
$e_0$ is a primitive element of the generalized  parabolic algebra
$\got{b}$  (cf. Definition IX.1.1 p. 328 in
\cite{neeb}).

Let us denote ${\mb{X}}:=dT(X), X\in \Ugc$, where    $\U$  denotes the
universal enveloping algebra. Let 
$\tilde{g}(g)=\exp X, \tilde{g}(g)\in G/H,~ X\in\got{m}$,
$e_{\tilde{g}(g)}=\exp({\mb{X}})e_0, X\in\got{m}$.
Let us remember again   
 Theorem XV.1.1 p. 646 in \cite{neeb}.
Note  that $T_o(G/H)\cong\got{g}/\got{h}\cong \gc/\bar{\got{b}}\cong
(\got{b}+\bar{\got{b}})/\bar{\got{b}}\cong
\got{b}/h_{\C}$,
 where we
have a linear isomorphism $\alpha:\got{g}/\got{h}\cong
\gc/\bar{\got{b}}$, $\alpha(X+\got{h})=X+\bar{\got{b}}$ (\cite{neeb1}). 
We can take
instead of $\got{m}\subset\g$ the subspace $\got{m}' \subset \gc$
complementary to $\bar{\bb}$, or the subspace of $\got{b}$
complementary to $\got{h}_{\C}$. Then let
 $x(s)$, $s\in [0,1]$,   be the one-parameter subgroup generated by
$X\in\m '$ and  
$x^*(s)$ his image in the reductive homogeneous space $M\cong 
G/H\hookrightarrow G_{\C}/B$, i.e. $x^*(s)=\exp (sX) .o$. 
If we choose a local {\em canonical} system of coordinates
$\{z_{\alpha}\}$ with respect to the basis $\{X_{\alpha}\}$ in \m',
then we can introduce the vectors
\begin{equation}\label{cvect}
e_{z}=\exp(\sum_{X_{\alpha}\in\got{m}'}z_{\alpha}{\mb{X}}_{\alpha}).e_0
\in\Hi .
\end{equation}
We get
\begin{equation}\label{unu2}
e_{\sigma (z)}=T(\sigma (z)), ~z\in M, 
\end{equation}
 and  we prefer to choose local coordinates such that
\begin{equation}\label{doi}
e_{\sigma (z)}=N(z)e_{\bar{z}},~~
N(z)=(e_{\bar{z}},e_{\bar{z}})^{-1/2}.
\end{equation}
 
Equations (\ref{cvect}), (\ref{unu2}), and (\ref{doi})
 define locally  the {\em coherent vector
 mapping}
\begin{equation}\label{cvm}
\varphi : M\rightarrow \bar{\Hi}, ~ \varphi(z)=e_{\bar{z}},  
\end{equation}
where $ \bar{\Hi}$ denotes the Hilbert space conjugate to $\Hi$.
We call the  vectors $e_{\bar{z}}\in\bar{\Hi}$ indexed by the points
 $z \in M $  {\it
Perelomov's coherent state vectors} \cite{perG}. Below we give more
details
about  
 this construction.

 Let  $V$ be a generalized highest weight module 
 with highest
weight $\lambda$ and primitive element $v_{\lambda}$.
 Then $V$ carries a non-degenerate
contravariant hermitian form if and only if $ V\cong\
L(\lambda , \bb):=
 M(\lambda ,\bb )/R $, where $M(\lambda ,\bb )$ is
the Verma module and $R$ is  the radical of the contravariant hermitian form
(cf. Definition IX.1.8 p. 333 in \cite{neeb}).
Now we suppose that $V$ is a highest weight module with respect to a
positive system of roots $\Delta^+$.
  Let
$\got{n}^{\pm}:=\sum_{\alpha\in\Delta^+}\g^{\pm\alpha}$ and
$e_0:=v_{\lambda}$ a primitive element. Let us now take into account
that
 $\Ug =\U
(\got{n}^-)\Ub$ and let us choose a canonical system of coordinates in the
highest weight module with respect to a fixed base of \g .
 One has locally finite representations by
direct exponentiation of the module $L(\lambda , \bb)$ (cf. Corollary
XII.2.7 p. 523  to the globalization  
Theorem XII.2.6  p. 521 in \cite{neeb}), and the
{\em Perelomov's coherent state vectors can be obtained by just taking
 the exponential of
images by the highest weight representation
of elements  of  $\U (\got{n}^-)$}.
 We can apply a {\em Gauss type decomposition}
as furnished by the Lemma XII.1.2. p. 495 in \cite{neeb}
 and get the coherent vector
mapping given by equations (\ref{unu2}), (\ref{doi}),
 (\ref{cvm}). 

\section{ The \FSB ~  \fl~ as reproducing kernel Hilbert space }\label{fock}

We have considered homogeneous CS-manifolds $M\cong G/H$ whose complex
structure comes from the embedding $i_1:M\hookrightarrow
G_{\C}/B$. We have chosen a section $\sigma: G_{\C}/B\to  G_{\C}$,
 and $G_{\C}$ can be regarded
as a complex analytic principal bundle     
$ B\stackrel{i}{\rightarrow}G_{\C}\stackrel{\lambda}{\rightarrow}
G_{\C}/B$.

 Let us introduce the function $f'_{\psi}:G_{\C}\to \C$
$f'_{\psi}(g):=(e_g,\psi ), g\in G, \psi\in \Hi$.
Then 
$f'_{\psi}(gb)=\chi(b)^{-1}f'_{\psi}(g), g\in G_{\C}, b\in B$,
where  $\chi$ is  the continuous homomorphism of the 
isotropy subgroup $B$ of  $G_{\C}$  
  in $\C^*$.   If    the homomorphism $\chi$  is
holomorphic, then
 {\it the coherent states realize the space
of holomorphic global sections} $
\Gamma^{\text{{hol}}}(M,L_{\chi})=H^0(M,L_{\chi})$ {\it on the}
 $G_{\C}$-{\it homogeneous line bundle
$L_{\chi}$ associated by means of the character
 $\chi$ to the  principal B-bundle}
 (cf. \cite{onofri}).
 Here the  holomorphic line bundle is $L_{\chi}:=M\times_{\chi}\C$,  
 also denoted   $L:=M\times_B\C$ (cf. \cite{bott,tirao}).

The local trivialization of the line bundle $L_{\chi}$ associates to
every $\psi\in\Hi$ a holomorphic function $f_{\psi}$   on a  open
set  in $M \hookrightarrow G_{\C}/B$. Let the notation $G_S:=G_{\C}
\setminus S$,
 where $S$ is the   set
$S:=\{g\in G_{\C}|\alpha_g=0\}$,
and $\alpha_g:= (e_g,e_0)$. $G_S$ is a dense subset of $G_{\C}$.
We  introduce  the function $f_{\psi}:G_{S}\mapsto\C$
$f_{\psi}(g)=\frac{f'_{\psi}(g)}{\alpha_g}, \psi\in\Hi ,~g\in G_{S}$.
The
function $f_{\psi}(g)$ on $G_S$
 is actually a function of the projection
$\lambda (g)$,
 holomorphic in  $M_S:=\lambda ( G_S)$.
We have introduced the function 
$f_{\psi}(g)=f_{\psi}(z_g)=\frac{(e_{\bar{z}_{g}},\psi )}
{(e_{\bar{z}_{g}},e_0)}
~~((e_{z_{g}},e_0)\not= 0)$
and also the coherent state map
$\varphi :M\rightarrow\overline{\Hi}^{\infty} , \varphi (z)=e_{\bar{z}},
z\in{\mc{V}}_0,$
where the canonical coordinates
$z=(z_1,\ldots ,z_n )$ constitutes
 a   local chart   on  ${\mathcal{V}}_0:=M_S\rightarrow
\C^n$,
such that $0=(0,\ldots ,0)$  corresponds to $\{ B\}$.
Note  also that
 ${\mathcal{V}}_0\equiv M\setminus \Sigma_0$, where $\Sigma_0:=\lambda
 (S)$ is the set of points of $M$ for which the
coherent state vectors are orthogonal to $e_0\in \Hi$, called 
{\it
polar divisor} of the point $z=0$  (cf. \cite{sbcl}).

Supposing that
{\it the line bundle $L_{\chi}$ is already very ample},
\fl~ is defined as 
$\fl :=\{f\in L^{2}(M,L)\cap{\mc{O}}(M,L)
|(f,f)_{\fl}<\infty \}$
 with respect to  the scalar product
\begin{equation}\label{scf}
(f,g)_{\fl} =\int_{M}\bar{f}(z)g(z)d\nu_M(z,\bar{z}),
\end{equation} 
where  $d{\nu}_{M}(z,\bar{z})$ is the invariant measure 
$\frac{d{\mu}_M(z,\bar{z})}{(e_{\bar{z}},e_{\bar{z}})}$, and
 $d{\mu}_M(z,\bar{z})$ represents the Haar measure on $M$.
It can be shown that the $\fl := L^{2,{\text{hol}}}(M,L_{\chi})$ {\em  is a closed
subspace of $L^2(M,L_{\chi})$ with continuous point evaluation}
(cf. \cite{PWJ}). 
 
Note that eq.  (\ref{scf})
is nothing else than the Parseval  ({\em overcompletness})
 identity  \cite{berezin}:
\begin{equation}\label{orthogk}
(\psi_1,\psi_2)=\int_{M=G/K}(\psi_1,e_{\bar{z}})(e_{\bar{z}},\psi_2)
d{\nu}_{M}(z,\bar{z}),~ (\psi_1,\psi_2\in \Hi ) .
\end{equation}

Let us  now introduce the map
\begin{equation}\label{aa}
\Phi :\Hi^{\star}\rightarrow \fl ,\Phi(\psi):=f_{\psi},
f_{\psi}(z)=\Phi(\psi )(z)=(\varphi (z),\psi)_{\Hi}=(e_{\bar{z}},\psi)_{\Hi},~
z\in{\mathcal{V}}_0,
\end{equation}
where we have identified the space  $\overline{\Hi}$ with the dual
space
$\Hi^{\star}$ of $\Hi$.

In fact, our supposition that $L_{\chi}$ is already a very ample line
bundle implies the validity of eq. (\ref{orthogk}) (cf. 
 Theorem XII.5.6 p. 542 in \cite{neeb},
 Remark VIII.5 in \cite{neeb94}, and 
 Theorem XII.5.14 p. 552 in \cite{neeb}).
 Rosenberg and Vergne  have shown
that the projectively induced 
 line bundle $\LL\cong L_{\chi}$ is ample, i.e. there exists
$n\in\N_0$ such that $L^{2,\text{hol}}(M,\LL^n )\not= \{0\}$ (cf.
Theorem 2.15 in  \cite{rv}; see also \S 4 in \cite{lis}).
  If the line
bundle is only a ample one, not every highest weight representation
leads to square integrable representations, and the highest weight
vector $e_0:=e_{\lambda}$ in the definition of coherent state vectors
has to verify a condition which generalizes the Harish-Chandra
condition in the semisimple case  (cf.
 Theorem XII.5.14 p. 552 in \cite{neeb} and 
 Remark VIII.5 in \cite{neeb94}).

It can be seen 
that  the group-theoretic relation (\ref{orthogk}) on homogeneous
manifolds fits into 
 Rawnsley's (global) realization \cite{raw} of Berezin's coherent states on 
quantizable K\"ahler manifolds \cite{berezin}. We emphasize that,
strictly speaking,    
{\it equation (\ref{orthogk})  should be
considered with a partition of unity.}

Let us introduce the notation
\begin{equation}\label{kernel1}
K:=\Phi\circ \varphi,~ K:M\rightarrow \fl ,~~~ 
%\end{equation}
%and denote
%\begin{equation}\label{kernel2}
K_w:=f_{e_{\bar{w}}}\in\fl.
\end{equation}
It can be defined a function, also denoted $K$,
   $K: M\times\overline{M}\rightarrow \C$, which on  ${\mathcal{V}}_0\times
\overline{\mathcal{V}}_0$ reads
\begin{equation}\label{kernel}
K(z,\overline{w}):=K_w(z)=
(e_{\bar{z}},e_{\bar{w}})_{\Hi}.
\end{equation} 
  For fixed $z\in M$  the function
 $K(z,\overline{w})$ is defined  for $w\notin \Sigma_z$
 \cite{sblb}. In the compact case $K(z,\overline{w})=0$ for $z\in M,~
 w\in\Sigma_z$.

Taking into account (\ref{aa}) and supposing that
eq. (\ref{orthogk})
is true, it follows that  the function $K$ (\ref{kernel}) is a reproducing
kernel. 
Using  the terminology of ref. \cite{neeb}, we have:
\newpage
\begin{Proposition}\label{realization}
Let $(T,\Hi)$ be a CS-representation and let us consider
 the Perelomov's CS-vectors defined in 
(\ref{cvect})-(\ref{doi}). Suppose that the  line bundle $L$ is
 very ample.  Then

\mbox{\rm{i)}} The function $K:M\times\overline{M}\rightarrow \C$, 
 $K(z,\overline{w})$ defined by equation   
(\ref{kernel}), is a   reproducing kernel.

\mbox{\rm{ii)}}  Let \fl~ be
the space $L^{2,\text{\em{hol}}}(M, L)$ endowed with the scalar
 product (\ref{scf}). Then 
 \fl~ is the reproducing kernel Hilbert space $\Hi_K\subset \C^M$
  associated to the kernel $K$
(\ref{kernel}).

\mbox{\rm{iii)}} The evaluation  map $\Phi$ defined in
 eqs. (\ref{aa}) 
extends to an  isometry   
\begin{equation}\label{anti}
(\psi_1,\psi_2)_{\Hi^{\star}}=(\Phi (\psi_1),\Phi
(\psi_2))_{\fl}=(f_{\psi_{1}},f_{\psi_{2}})_{\fl}=
\int_M\overline{f}_{\psi_1} (z)f_{\psi_2}(z)d\nu_M(z),
\end{equation}
and  the overcompletness eq. (\ref{orthogk}) is verified. 
\end{Proposition}

\section{Representations of CS-Lie algebras  by
 differential operators}\label{DIFF}

We remember the definitions of the functions $f'_{\psi}$  and
$f_{\psi}$,  which allow to
write down
\begin{equation}\label{iar}
f_{\psi}(z)=(e_{\bar{z}},\psi)=\frac{(T(\bar{g})e_0,\psi)}
{(T(\bar{g})e_0,e_0)},~ z\in M,~ \psi\in \Hi .
\end{equation}
So, we get 
\begin{equation}\label{iar1}
 f_{T(\overline{g'}).\psi}(z)= \mu
(g',z)f_{\psi}(\overline{g'}^{-1}.z),
\end{equation}
where 
\begin{equation}\label{iar2}
 \mu (g',z)=
\frac{(T(\overline{g'}^{-1}\overline{g})e_0,e_0)}{(T(\overline{g})e_0,e_0)}
=\frac{\Lambda (g'^{-1}g)}{\Lambda (g)}.
\end{equation}
We remember that 
${T(g).e_0=e^{i\alpha (g)}e_{\tilde{g}}=\Lambda (g) e_{z_{g}}}$
where we have used the decompositions
$g=\tilde{g}.h, ~(G=G/H.H);~~ g = z_g. b ~(G_{\C}=G_{\C}/B.B)$.
 We have also  the relation
$\chi_0(h)=e^{i\alpha (h)},~ h\in H$ and 
$ \chi (b) =\Lambda (b),~b\in B$, where
$\Lambda(g)=\frac{e^{i\alpha(g)}}{(e_{\bar{z}},e_{\bar{z}})^{1/2}}$.
We can also write down another expression for multiplicative factor
 $\mu$ appearing  in eq. (\ref{iar1})    using the CS-vectors
\begin{equation}\label{iar3}
\mu (g',z)=\Lambda(\bar{g'})(e_{\bar{z}},e_{\bar{z'}})=e^{i\alpha
  (\bar{g'})}
\frac{(e_{\bar{z}},e_{\bar{z'}})}{(e_{\bar{z'}},e_{\bar{z'}})^{1/2}}.
\end{equation}

The following assertion is easy to be checked using successively
eq. (\ref{iar2}):
\begin{Remark}Let us consider the relation (\ref{iar}). Then we have
(\ref{iar1}), where  $\mu$ can be written down as in equations
(\ref{iar2}), 
(\ref{iar3}). 
We have the relation
 $\mu (g,z) =J(g^{-1},z)^{-1}$, i.e. the multiplier $\mu$
 is the cocycle in the  unitary representation
$(T_K,\Hi_K)$ attached to the
positive definite holomorphic kernel $K$ defined by equation (\ref{kernel}),
 \begin{equation}\label{num}
(T_K(g).f)(x):=J(g^{-1},x)^{-1}.f(g^{-1}.x),
\end{equation}
and the cocycle verifies the relation
\begin{equation}\label{prod}
J(g_1g_2,z)=J(g_1,g_2z)J(g_2,z).
\end{equation}
\end{Remark}
  Note that
{\em   the prescription (\ref{num})
 defines
a continuous action of $G$ on ${\mbox{\rm{Hol}}}(M,\C )$ with respect to
the compact open topology on the space  ${\mbox{\rm{Hol}}}(M,\C )$.
 If $K:M\times M\rightarrow \C$ is a continuous positive definite kernel
holomorphic in the first argument satisfying
$K(g.x,\overline{g.y})=J(g,x)K(x,\overline{y})J(g,y)^*$,
$g\in G$, $x,y\in M$, then the action of $G$ leaves the reproducing kernel
Hilbert space $\Hi_K\subseteq {\mbox{\rm{Hol}}}(M,\C )$ invariant and defines
a continuous unitary representation $(T_K,\Hi_K)$ on this space}
(cf. Prop. IV.1.9 p. 104 in  Ref. \cite{neeb}).

Let us consider the triplet $(G, T, \Hi )$. Let
 $\Hi^0:=\Hinf$, considered as a pre-Hilbert space, and
let  $B_0(\Hi^0)\subset \mc{L}(\Hi )$ denote the set
  of linear operators
 $A:\Hi^0\rightarrow \Hi^0$
which have  a formal adjoint  $A^{\sharp}:\Hi^0\rightarrow\Hi^0$, i.e.   
$(x,Ay)=(A^{\sharp}x,y)$ for all $x,y\in \Hi^0$.
Note that if  $B_0(\Hi^0)$ is the set of unbounded operators on \Hi ,
then the domain $\mathcal{D}(A^*)$ contains $\Hi^0$ and $A^*\Hi^0\subseteq 
 \Hi^0$, and  it make sense to refer to the closure $\overline{A}$
 of $A\in B_0(\Hi^0)$ (cf. \cite{neeb} p. 29; here $A^*$ is the
 adjoint of $A$).

Let \g~ be the Lie algebra of $G$ and let us denote
 by $\mathcal{S}:=\Ugc$ the semigroup associated with
 the universal enveloping  algebra equipped with the antilinear involution
 extending the antiautomorphism $X\mapsto X^*:=-\bar{X}$ of $\gc$.  {\it
 The derived representation} is defined as
\begin{equation}\label{derived}
dT:\Ugc\rightarrow B_0(\Hi^0),~~ \mbox{\rm{with}} ~~ dT(X).v
:=\left.\frac{d}{dt}\right|_{t=0}T(\exp tX).v, X\in\g .
\end{equation}
Then $dT$ is a {\it hermitian representation}
 of $\mathcal{S}$ on $\Hi^0$ (cf. Neeb \cite{neeb}, p. 30).
Let us denote his image in $B_0(\Hi^0)$ with $\am := dT(\mathcal{S})$. 
If $\Phi : \Hi^{\star}\rightarrow \fl $ is the (Segal-Bargmann) isometry
(\ref{aa}), we are interested in the study of the image of 
\am~  via $\Phi$  as subset in the
algebra of holomorphic, linear differential operators,  
$ \Phi\am\Phi^{-1}:={\db{A}}_M\subset\got{D}_M$.
The new results for semisimple Lie groups $G$ for
 ${\db{A}}_M$ , $M\approx G/H$,
 are contained in the main theorem in \cite{sbcpol,sbctim}.

The {\it  sheaf}
 $\got{D}_M$ (or simply \D ) {\it of holomorphic, finite order, linear
differential operators on} $M$ is a
 subalgebra of homomorphism ${\mathcal Hom}_{\C}({\cal O}_M,{\cal O}_M)$ 
 generated
 by the sheaf ${\cal O}_M$ of germs of holomorphic functions of $M$ and the
 vector fields. 
 We consider also {\it  the subalgebra} \AM~ of ${\db{A}}_M$~
 {\it of differential operators with
 holomorphic polynomial coefficients}.
Let $U:=\mathcal{V}_0$ in $M$, endowed with the 
coordinates
$(z_1,z_2,\cdots ,z_n)$. We set $\pa_i:=\frac{\pa}{\pa z_i}$ and
$\pa^{\alpha}:=\pa^{\alpha_1}_1 
\pa^{\alpha_2}_2\cdots \pa^{\alpha_n}_n$, $\alpha :=(\alpha_1,
\alpha_2 ,\cdots ,\alpha_n)\in\N^n$. The sections of \DM~ on $U$ are
$A:f\mapsto \sum_{\alpha}a_{\alpha}\pa^{\alpha}f$,
$a_{\alpha}\in\Gamma (U, {\cal{O}})$, the $a_{\alpha}$-s being zero
except a finite number. 

For $k\in\N$, let us denote by $\D_k$ the subsheaf of differential
operators of degree $\le k$ and by $\D'_k$ the subsheaf of elements of 
$\D_k$ without constant terms.  $\D_0$ is identified with
$\cal{O}$ and $\D'_1$ with the sheaf of vector fields. The
filtration of \DM~ induces a filtration on  ${\A}_M$.

Summarizing, we have the following three objects which correspond
each to other:
\begin{equation}\label{correspond}
\g \ni X \mapsto\mb{X}\in\am\mapsto\db{X}\in\AAA\subset \DM, {\text
  {~differential  operator  on}}~ \fl .
\end{equation}
Now we can see that
 \begin{Proposition}If $\Phi$ is the isometry
(\ref{aa}), then 
$\Phi dT(\g_{\C})\Phi^{-1}\subseteq \D_1$.
\end{Proposition}
{\it Proof}. Let us consider an element in $\g_{\C}$ and his image in \DM ,
via the correspondence (\ref{correspond}):
$$\g_{\C}\ni G \mapsto{\db{G}}\in \DM;~~
 {\db{G}}_z(f_{\psi}(z))= 
\db{G}_z(e_{\bar{z}},\psi)= (e_{\bar{z}},\mb{G}\psi),$$
$$\mb{G}=dT (G)=\frac{d}{dt }|_{t=0}T(\exp (tG)).$$
Remembering equation
(\ref{iar1}) and  determining  the derived representation,
 we get
\begin{equation}\label{sss}
\db{G}_{z}(f_{\psi}(z))=(P_{G}(z)+
\sum Q^i_G(z)\frac{\pa}{\pa  z_i})f_{\psi}(z);~
\end{equation}
 $$P_G(z)=\frac{d}{dt}|_{t=0}\mu (\exp (tG),z);~ 
 Q^i_G(z)=\frac{d}{dt}|_{t=o}(\exp(-tG).z)_i.~~~~~~~~~~~~~ \gata~~~$$
Now we formulate the following assertion:
\begin{Remark}\label{main}
If $(G,T)$ is a CS-representation, then \AAA~ is a subalgebra of
holomorphic differential operators with polynomial coefficients,
 $\AAA \subset \AM\subset \DM$. 
%\enlargethispage{1cm}
More exactly, for    $X\in\g$,  
  let us denote by $\mb{X}:= dT(X)\in\am $, where the action is
considered on the space of functions \fl . Then, for 
CS-representations, $\db{X}\in \A_1=\A_0\oplus \A_1'$.
 
Explicitly, if  $\lambda\in\Delta $ is a root and
 $G_{\lambda}$ is in a base of the Lie algebra $\g_{\C}$ of $G_{\C}$, then
 his image  ${\db G}_{\lambda}\in\DM$ acts as a first order
 differential operator on the
\FSB \fl
  \begin{equation}\label{generic}
{\db G}_{\lambda}= P_{\lambda}+\sum_{\beta \in\Delta_{\m '}} Q_{\lambda
 ,\beta} \pa_{\beta},~ \lambda \in \Delta,
 \end{equation}
where $ P_{\lambda}$ and $ Q_{\lambda ,\beta}$  are 
 polynomials in
 $z$, and $\m '$ is the subset of $\g_{\C}$ which appears in the
definition  (\ref{cvect}) of the coherent state vectors.
\end{Remark}
 Actually, we don't have a proof of this 
assertion for the general case of CS-groups.
  For the compact case, there exists the calculation of
Dobaczewski \cite{dob}, which in fact can be extended
also to real semisimple Lie algebras. For compact hermitian symmetric
spaces it was shown \cite{sbcag} that degrees of the polynomials
 $P$ and $Q$-s are   
 $\le 2$ and similarly for the non-compact hermitian symmetric  case
\cite{sbl}.
Neeb \cite{neeb} gives a proof of this Remark
  for 
CS-representations for  the (unimodular) Harish-Chandra type
groups. Let us also remember that:
{\it If $G$ is an admissible Lie group such that the universal
complexification $G\rightarrow G_{\C}$ is injective and $ G_{\C}$ is
simply connected, then $G$ is of Harish-Chandra type}
(cf. Proposition V.3 in \cite{neeb94}). The derived
 representation (\ref{derived}) is obtained
   differentiating eq. (\ref{num}), and we get two 
 terms, one in $\D_0$ and the
other one in $\D_1'$.
A proof that the two parts are in fact $\A_0$ and respectively $\A_1'$
    is contained  in
Prop. XII.2.1 p. 515 in \cite{neeb}
  for the groups of Harish-Chandra type in the particular situation
where the space $\got{p}^+$ in  Lemma VII.2.16 p. 241 in \cite{neeb}
 is abelian. We have presented 
 explicit formulas for semisimple  Lie groups and also the 
simplest example where the maximum  degree of the polynomials 
$P$ and $Q$-s is  3 (cf. \cite{sbcpol,sbctim}).
 \hfill   $\gata$

{\bf Acknowledgment} S.B. is grateful to Karl-Hermann Neeb  for many
suggestions, criticism and the proof of Remark \ref{red}. This
investigation was partially supported by the MEC project CERES 66/2001.

\end{document}